\documentclass[bibay2,jsl]{asl}
\usepackage{amssymb,
  bussproofs,
  graphicx,
  bussproofs-extra,
  stmaryrd}
\let\citet\cite
\let\citep\cite
\let\citealt\cite
\let\citeauthor\citeauth

\usepackage[hidelinks]{hyperref}

\makeatletter
\def\Rightarrowfill@#1{\m@th\setboxz@h{$#1=$}\ht\z@\z@
  $#1\copy\z@\mkern-6mu\cleaders
  \hbox{$#1\mkern-2mu\box\z@\mkern-2mu$}\hfill
  \mkern-6mu\mathord\Rightarrow$}
\newcommand{\xRightarrow}[2][]{%
  \mathrel{\mathop{%
    \setbox\z@\vbox{\m@th
      \hbox{$\scriptstyle\;{#1}\;\;$}%
      \hbox{$\m@th\scriptstyle\;{#2}\;\;$}%
    }%
    \hbox to\ifdim\wd\z@>\minaw@\wd\z@\else\minaw@\fi{%
      \Rightarrowfill@\displaystyle}}%
  \limits^{#2}\@ifnotempty{#1}{_{#1}}}%
}\makeatother

\newtheorem{thm}{Theorem}

\newtheorem{prop}[thm]{Proposition}

\newtheorem{cor}[thm]{Corollary}
\theoremstyle{definition}
\newtheorem{definition}[thm]{Definition}
\newtheorem{example}[thm]{Example}

\def\HSeq{\ensuremath{\mathrel{\,\sslash\,}}}
\def\Red{\mathrm{Red}}

\def\Succ{\mathrm{Succ}}

\let\models\vDash
\let\nmodels\nvDash
\newcommand{\Seq}[1][]{\ensuremath{\mathrel{\xRightarrow{#1}}}}
\let\phi\varphi

\RequirePackage{tikz}
\usetikzlibrary{arrows,automata,intersections,positioning}
\tikzset{
  modal/.style={>=stealth',
    shorten >=1pt,
    shorten <=1pt,
    auto,
    node distance=1.5cm,
    label distance=2pt,
    semithick},
  every label/.style={phantom,align=left},
  world/.style = {circle,draw,minimum size=0.5cm,fill=gray!15},
  modal every node/.style={world},
  point/.style={circle,draw,inner sep=0.5mm,fill=black},
  phantom/.style={rectangle,inner sep=0pt,draw=none,fill=none},
  reflexive above/.style={->,loop,looseness=7,in=60,out=120},
  reflexive below/.style={->,loop,looseness=7,in=240,out=300},
  reflexive left/.style={->,loop,looseness=7,in=150,out=210},
  reflexive right/.style={->,loop,looseness=7,in=30,out=330}
}

\title[Modular Hypersequent Calculi for K, T, and D]{Cut-free Completeness for Modular Hypersequent Calculi for Modal Logics K, T, and D}

\author{Samara Burns}
\revauthor{Burns, Samara}
\address{Columbia University\\
Department of Philosophy\\
1150 Amsterdam Avenue \\
New York, NY 10027, USA}
\email{sb4318@columbia.edu}

\author{Richard Zach}
\revauthor{Zach, Richard}
\address{University of Calgary\\
Department of Philosophy\\
2500 University Drive NW\\
Calgary, AB T2N 1N4, Canada}
\email{rzach@ucalgary.ca}
\urladdr{https://richardzach.org/}
\thanks{Forthcoming in \emph{The Review of Symbolic
Logic},
\href{https://doi.org/10.1017/S1755020320000180}{DOI:10.1017/S1755020320000180}.
\copyright~Cambridge University Press.}
\begin{document}

\def\journalname{The Review of Symbolic Logic}
\def\shortjournalname{The Review of Symbolic Logic}
\def\longjournalname{The \hfil Review \hfil of \hfil Symbolic \hfil
Logic}
\def\jslname{{\bfseries\itshape\selectfont The Journal of Symbolic Logic}}
\def\ISSN{1755-0203}

\begin{abstract}
We investigate a recent proposal for modal hypersequent calculi. The
interpretation of relational hypersequents incorporates an
accessibility relation along the hypersequent. These systems give the
same interpretation of hypersequents as Lellman's linear nested
sequents, but were developed independently by Restall for S5 and
extended to other normal modal logics by Parisi. The resulting systems
obey Do\v{s}en's principle: the modal rules are the same across
different modal logics. Different modal systems only differ in the
presence or absence of external structural rules. With the exception
of S5, the systems are modular in the sense that different structural
rules capture different properties of the accessibility relation. We
provide the first direct semantical cut-free completeness proofs for
K, T, and D, and show how this method fails in the case of B and S4.
\end{abstract}

\maketitle 

\section{Introduction}

Modal sequent calculi have been developed for~K and many of its
extensions, but it has historically been difficult to develop sequent
systems that have nice proof-theoretic properties. The cut elimination
theorem is an important result in structural proof theory: any sequent
that is derivable in a calculus can be derived without the use of cut.
Notably, the sequent system for~S5 given by \citet{Ohnishi1957} is not cut-free.
Although cut-free sequent systems for S5 were later developed by
\citet{Ohnishi1982} and \citet{Brauner2000}, this issue prompted
research into extensions of the sequent calculus that could better
accommodate modal logics. One such extension are hypersequent calculi,
which operate on sets or sequences of sequents. The first hypersequent
system, also for~S5, was presented by \citet{Mints1971,Mints1974}.
There has since been a proliferation of hypersequent approaches to
modal logics (\citealt{Pottinger1983, Avron1996, Brunnler2009,
Restall2009, Poggiolesi2008, Indrzejczak2012,
Lahav2013}).

Sequents $\Gamma \Seq \Delta$ can be translated into single formulas:
$\bigwedge \Gamma \to \bigvee \Delta$, or sometimes $\Box(\bigwedge
\Gamma \to \bigvee \Delta)$ in the case of sequent systems for modal
logics. Earlier hypersequent approaches to modal logics interpret
hypersequents as disjunctions of the formula interpretations of the
individual sequents. Intuitively, a hypersequent is evaluated at a
single world, and describes a disjunction. Under this interpretation,
the order and multiplicity of sequents in a hypersequent is
immaterial, and external contraction, exchange, and weakening are
admissible rules. More recent approaches, however, interpret
hypersequents in such a way that different sequents in a hypersequent
are evaluated at different worlds. If the logic is S5 (and so either
no or a universal accessibility relation is assumed), the order and
multiplicity of sequents in a hypersequent still does not matter. But
for other logics, the worlds at which adjacent sequents in a
hypersequent are evaluated must be related. Then the order and
multiplicity of sequents in a hypersequent \emph{does} matter, and
external structural rules such as external contraction, weakening, and
exchange are not sound in general. Approaches using this interpretation
of hypersequents are the linear nested
sequent systems of \citet{Lellmann2015}, \citet{Lellmann2016}, and
\citet{GoreLellmann2019}, the non-commutative hypersequents of
\citet{Indrzejczak2016, Indrzejczak2018, Indrzejczak2019}, and the
ordered hypersequents of
\citet{BaeldeLickSchmitz2018a}. The 2-sequents of \citet{Masini1992}
are a notational variant of hypersequents with the same interpretation
as linear nested sequents. \citet{KuznetsLellmann2018} applied the
linear nested sequent approach also to G\"odel logic. 

The hypersequent system for S5 of \citet{Restall2009}, although taking
hypersequents as sets of sequents, explicitly interprets individual
sequents as describing different possible worlds. \citet{Parisi2017}
generalized this interpretation by incorporating an accessibility
relation into the interpretation of a hypersequent, and offered
calculi for K, T, D, S4, and~S5. The interpretation of Parisi's
hypersequents is equivalent to the formula interpretation of linear
nested sequents. To unify terminology, we'll call hypersequents
\emph{relational} if their interpretation takes the accessibility
relation into account.

In addition to cut elimination, there are other desiderata that one
might consider when developing hypersequent calculi.  We consider two
properties of hypersequent systems that have been proposed as
important. One is \emph{modularity:} each property of the
accessibility relation is captured by a single rule or set of rules.
Modularity yields a satisfying systematicity for proof systems for
various kinds of modal logics. The fact familiar from modal
correspondence theory that properties of the accessibility relation
can be captured by different modal axioms allows for elegant treatment
of large classes of logics and uniform results. Likewise, analogous
modularity of proof systems opens up the possibility of dealing with
combinations of properties of the accessibility relation not
piecemeal, but systematically by combining different structural rules.
Another property relational hypersequent systems have is that they
obey what's been called Do\v{s}en's principle: hypersequent systems
for different modal logics only differ in the presence or absence of
structural rules, while the rules for modalities are the same for all
systems. This corresponds to a methodological principle that the
meaning of a connective should be determined by its rules of
inference. So, the rules for modal operators should be the same
regardless of the structure of the accessibility relation.

Parisi's systems are the first candidates for hypersequent calculi for
modal logics that are both modular and conform to Do\v{s}en's
principle. These systems are unique in that they do not require the
addition of rules that govern the modal operator when moving between
modal systems. Different modalities are treated using only external
structural rules.\footnote{The system of \citet{Masini1992} for~KD
uses an implicit Drop rule. \citet{Lellmann2015} proposes the same
structural rules (drop, external weakening) to deal with D and~4, but
does not study the resulting calculi in detail.} This is not to say
that modular systems obeying Do\v{s}en's principle are superior to
other approaches. In fact, e.g., \cite[32--34]{Poggiolesi2010} has
called Do\v{s}en's principle into question, and perhaps the ideal of
modularity simply cannot be universally combined with other results
such as cut elimination.

The prospects of the relational hypersequent approach to obtain
cut-free systems that are modular and obey Do\v{s}en's principle is
nevertheless an important and interesting question, which this paper
aims to shed light on. In section~\ref{sec:sound}, we introduce
Parisi's relational hypersequent approach and compare it to Lellman's
linear nested sequents. \citet{Lellmann2015} and \citet{Parisi2017}
showed completeness for their systems by inductively translating usual
sequent derivations into hypersequent derivations. Cut-free
completeness for the hypersequent systems then follows from the
cut-free completeness of the corresponding sequent system. Exceptions
are Parisi's systems RB, RS4, and RS5, where the translation makes use
of the cut rule. \citet{Masini1992} and \citet{Parisi2017}
independently gave syntactic cut-elimination proofs for their
(equivalent) systems for~D. \citet{Restall2009} shows cut-free
completeness for RS5 directly. We improve on these results by
providing a direct cut-free completeness proof for~RK, the relational
hypersequent calculus for~K (Section~\ref{sec:comp}). In
Section~\ref{sec:comp2} we show how this proof can be adapted to
obtain cut-free completeness for T and~D. Finally, in
Section~\ref{sec:fail} we discuss the limitations of the method for
the cases of logics B and~S4.

\section{Relational Hypersequent Calculi}\label{sec:sound}

\begin{definition}
We call any expression of the form $\Gamma \Seq \Delta$ a
\emph{sequent}, where $\Gamma$ and $\Delta$ are sets of formulas.

A \emph{hypersequent} is any expression of the form $S_1 \HSeq \dots
\HSeq S_n$, where the $S_i$ are sequents.
\end{definition}

The rules for the calculus RK, sound and complete for~K, are found in
table~\ref{table:RHS2}. To keep the subsequent proofs simple, we give
only the logical rules for $\lnot$, $\land$, and~$\Box$. (The rules for
$\lor$ and $\to$ are as usual, and rules for $\Diamond$ are
symmetrical to those for~$\Box$.) Below is an example proof in RK.
\begin{center}\def\fCenter{}
\Axiom$\fCenter\phantom{\HSeq {}} \phantom{{} \land \psi} \varphi  \Seq \varphi$ 
\RightLabel{EWL} 
\UnaryInf$\quad \Seq \phantom{\Box\varphi} \fCenter \HSeq \phantom{{}
\land \psi} \varphi
 \Seq \varphi$
\RightLabel{$\land$ L} 
\UnaryInf$\quad \Seq \phantom{\Box\varphi} \fCenter \HSeq \varphi \land \psi \Seq \varphi$
\RightLabel{$\Box$ L} 
\UnaryInf$\Box (\varphi \land \psi) \Seq \phantom{\Box\varphi}\fCenter \HSeq \phantom{\varphi \land \psi} \Seq
\varphi$
\RightLabel{$\Box$ R} 
\UnaryInf$\Box (\varphi \land \psi) \Seq \Box \varphi \fCenter$
\DisplayProof
\end{center}

Calculi for extensions of K are obtained by adding external structural
rules, which each characterize a property of the accessibility
relation. The structural rules and resulting calculi are summarized in
Tables \ref{table:RHS} and~\ref{table:RHS3}. As usual, we write
$\vdash_\mathrm{RX} H$ to mean that $H$ has a derivation in calculus~RX,
where X is one of K, T, 4, B, D, S4, S5.\footnote{Note
that we reverse the order of hypersequents in \cite{Parisi2017} to
facilitate comparison with the linear nested sequents of
\citet{Lellmann2015}.}

With the exception of RS5, these systems are modular: each external
structural rule represents an axiom characteristic of a property of
the accessibility relation. In the case of system RS5, the EE rule
does not only capture symmetry, but also transitivity. This way of
capturing S5 is equivalent to that of \citet{Restall2009}, but one may
also straightforwardly replace the EE rule with Sym. The resultant
calculus would be complete \citep{Parisi2017}, but the simulation of
sequent derivations in the hypersequent calculus uses cut.

\begin{table}[tp]
\caption{Rules of RK}
\label{table:RHS2}
\def\arraystretch{2.5}
\begin{tabular}{ll}
\hline\hline\vspace{-2ex}
\textbf{Axioms} &
$\phi \Seq \phi$\\
 
\textbf{Internal Structural Rules} &

\\

\rule{0pt}{4ex}   
\def\fCenter{\Seq}
\Axiom$G \HSeq \phantom{\varphi,} \Gamma \fCenter \Delta \HSeq H$
\RightLabel{WL}
\UnaryInf$G \HSeq \varphi, \Gamma \fCenter \Delta \HSeq H$
\DisplayProof &

\def\fCenter{\Seq}
\Axiom$G \HSeq \Gamma \fCenter \Delta \phantom{, \varphi} \HSeq H$
\RightLabel{WR}
\UnaryInf$G \HSeq \Gamma \fCenter \Delta, \varphi \HSeq H$
\DisplayProof 
\\


\multicolumn{2}{c}{
\def\fCenter{\Seq}
\Axiom$G \HSeq \Gamma \fCenter \Delta, \varphi \HSeq H$
\Axiom$G \HSeq \varphi, \Lambda \fCenter \Theta\HSeq H$
\RightLabel{Cut}
\BinaryInf$G \HSeq \Gamma, \Lambda \fCenter \Delta, \Theta \HSeq H$
\DisplayProof}

\\

\textbf{External Structural Rules} &

\\

\AxiomC{$G\phantom{{} \HSeq \Gamma \Seq \Delta}$}
\RightLabel{EWR}
\UnaryInfC{$G \HSeq \quad \Seq \quad$}
\DisplayProof &

\AxiomC{$\phantom{\Gamma \Seq \Delta \HSeq {}}G$}
\RightLabel{EWL}
\UnaryInfC{$\quad \Seq \quad \HSeq G$}
\DisplayProof
\\

\textbf{Logical Rules}&

\\

\def\fCenter{\Seq}
\Axiom$G \HSeq \phantom{\lnot \varphi,} \Gamma \fCenter \Delta,
\varphi \HSeq H$
\RightLabel{$\lnot$L}
\UnaryInf$G \HSeq \lnot\varphi, \Gamma \fCenter \Delta
\phantom{,\varphi} \HSeq H$
\DisplayProof &

\def\fCenter{\Seq}
\Axiom$G \HSeq \varphi, \Gamma \fCenter \Delta \phantom{,\lnot\varphi}
\HSeq H$
\RightLabel{$\lnot$R}
\UnaryInf$G\HSeq \phantom{\varphi,} \Gamma \fCenter \Delta,
\lnot\varphi\HSeq H$
\DisplayProof 
\\

\def\fCenter{\Seq}
\Axiom$G\HSeq \phantom{\land\psi,} \varphi, \Gamma \fCenter
\Delta\HSeq H$
\RightLabel{$\land$L$_1$}
\UnaryInf$G\HSeq \varphi \land \psi, \Gamma \fCenter \Delta \HSeq H$
\DisplayProof 

&

\def\fCenter{\Seq}
\Axiom$G\HSeq \phantom{\varphi\land,}\psi, \Gamma \fCenter
\Delta\HSeq H$
\RightLabel{$\land$L$_2$}
\UnaryInf$G\HSeq \varphi \land \psi, \Gamma \fCenter \Delta \HSeq H$
\DisplayProof \\

\multicolumn{2}{c}{
\def\fCenter{\Seq}
\Axiom$G\HSeq \Gamma \fCenter \Delta, \varphi\HSeq H$
\Axiom$G\HSeq \Gamma \fCenter \Delta, \psi\HSeq H$
\RightLabel{$\land$R}
\BinaryInf$G \HSeq \Gamma \fCenter \Delta, \varphi \land \psi \HSeq H$
\DisplayProof 
}

\\

\def\fCenter{\HSeq}
\Axiom$G \HSeq \phantom{\Box \varphi,} \Gamma \Seq \Delta
\fCenter \varphi, \Lambda \Seq \Theta \HSeq H$
\RightLabel{$\Box$L}
\UnaryInf$G \HSeq  \Box \varphi, \Gamma \Seq \Delta \fCenter
\phantom{\varphi,} \Lambda \Seq \Theta \HSeq H$
\DisplayProof&

\def\fCenter{\HSeq}
\Axiom$H \fCenter \Gamma \Seq \Delta \phantom{\Box \varphi,}
\HSeq \quad \Seq \varphi$
\RightLabel{$\Box$R}
\UnaryInf$H \fCenter \Gamma \Seq \Delta, \Box \varphi$
\DisplayProof 
\\
\hline\hline
\end{tabular}
\end{table} 

\begin{table}[t]
\caption{Additional external structural hypersequent rules}
\label{table:RHS}
\centering
\begin{tabular}{@{}c@{\qquad}l@{\qquad}l}
\hline\hline
Rule & Sound for & Axiom \\ \hline \\
\def\fCenter{\HSeq}
\AxiomC{$G \fCenter  \Gamma \Seq \Delta \HSeq \Gamma
  \Seq \Delta \HSeq H$}
\RightLabel{EC}
\UnaryInfC{$G \fCenter \Gamma \Seq \Delta \HSeq H$} \DisplayProof

& 
reflexive
&
$\Box A \Seq A$ (T)\\[3ex]

\AxiomC{$G \HSeq H$}
\RightLabel{EW}
\UnaryInfC{$G \HSeq \quad \Seq \quad \HSeq H$}
\DisplayProof \vspace{1mm}
&
transitive
& $\Box A \Seq \Box\Box A$ (4)
 \\[3ex]

\AxiomC{$G \HSeq  \quad \Seq \quad$}
\RightLabel{Drop} 
\UnaryInfC{$G \phantom{\HSeq  \quad \Seq \quad}$} 
\DisplayProof 

& serial & $\Box A \Seq \Diamond A$ (D)
\\[3ex] 

\AxiomC{$\Gamma_1 \Seq \Delta_1 \HSeq \ldots \HSeq \Gamma_n
\Seq \Delta_n$} \RightLabel{Sym} \UnaryInfC{$\Gamma_n \Seq
\Delta_n \HSeq \ldots \HSeq \Gamma_1 \Seq \Delta_1$}
\DisplayProof 
& symmetric

& $A \Seq \Box\Diamond A$ (B)
\\[3ex]

 \def\fCenter{\HSeq}
\Axiom$G \HSeq \Gamma \Seq \Delta
\fCenter \Pi \Seq \Lambda \HSeq H$
\RightLabel{EE}
\UnaryInf$G \HSeq \Pi \Seq \Lambda \fCenter
\Gamma \Seq \Delta \HSeq H$
\DisplayProof \vspace{1mm}

& symmetric and transitive\\
\hline\hline
\end{tabular}
\end{table}

\begin{table}[h]
\caption{Hypersequent calculi for various logics}
\label{table:RHS3}
\centering
\begin{tabular}{l@{\quad}l@{\quad}l}
\hline\hline
Calculus & Logic & External Structural Rules \\ \hline
RT & $T = KT$ & RK + EC \\
RB & $B = KB$ & RK + Sym \\
RD & $D = KD$ & RK + Drop \\ 
R4 & $4 = K4$ & RK + EW  \\
RS4 & $S4 = KT4$ & RK + EC + EW  \\ 
RS5 & $S5 = KT4B$ & RK + EC + EW + EE \\ 
\hline\hline
\end{tabular}
\end{table}

The semantics of relational hypersequents are given in terms of
absence of counterexamples.

\begin{definition}[Branch of worlds] Let $\mathfrak{F} = \langle W,
R\rangle$ be a frame. A \emph{branch of worlds} in $\mathfrak{F}$ is a
sequence $w_1$, \ldots,~$w_n$ of worlds such that $w_i R w_{i+1}$ for
$i=1$, \dots,~$n-1$.
\end{definition}

\begin{definition}[Countermodel]
A model $\mathfrak{M}$ is a \emph{countermodel} to a sequent $\Gamma
\Seq \Delta$ at a world $w$ iff for all $\varphi \in \Gamma,
\mathfrak{M}, w \models \varphi$ and for all $\psi \in \Delta,
\mathfrak{M}, w \nmodels \psi$.
\end{definition}

\begin{definition}[Counter-example] A model $\mathfrak{M}$ is a
\emph{counter-example} to a hypersequent $\Gamma_1 \Seq \Delta_1 \HSeq
\ldots \HSeq \Gamma_n \Seq \Delta_n$ iff there is a branch of worlds
$w_1$, \dots,~$w_n$ such that $\mathfrak{M}$ is a countermodel to
$\Gamma_i \Seq \Delta_i$ at $w_i$ for all $1 \leq i \leq n$.
\end{definition}

\begin{definition}[Valid hypersequent] A hypersequent $H$ is
\emph{valid} in a class of frames $\mathfrak{F}$ just in case there is
no counter-example to it that is in $\mathfrak{F}$. Otherwise, we say
that the hypersequent is \emph{invalid}.
\end{definition}

It is important to note that the interpretation of the relational
hypersequent is equivalent to that of the linear nested sequent, as in
\citet{Lellmann2015}. Whereas the interpretation of the relational
hypersequent is given in terms of a branch of worlds along a
hypersequent, linear nested sequents interpret the hypersequent as a
disjunction of nested modal formulas. This interpretation is given by
a mapping~$I$ on hypersequents~$H$:
\begin{align*}
I(\Gamma \Seq \Delta) & = \bigwedge \Gamma \to \bigvee \Delta\\
I(\Gamma \Seq \Delta \HSeq H) & = (\bigwedge \Gamma \to \bigvee \Delta) \lor \Box I(H)
\end{align*}

In other words, a linear nested sequent 
\[\Gamma_1 \Seq \Delta_1 \HSeq \dots \HSeq \Gamma_n \Seq \Delta_n\]
is interpreted as
\[
\bigwedge (\Gamma_1 \to \bigvee \Delta_1) \lor \Box((\bigwedge
\Gamma_2 \to \bigvee \Delta_2) \lor \Box(\ldots \Box (\bigwedge
\Gamma_n \to \bigvee \Delta_n) \ldots))
 \] 

\begin{prop}
A relational hypersequent $\Gamma_1
\Seq \Delta_1 \HSeq
\ldots \HSeq \Gamma_n \Seq \Delta_n$ is valid iff $I(\Gamma_1
\Seq \Delta_1 \HSeq \ldots \HSeq \Gamma_n \Seq \Delta_n)$
is valid.
\end{prop}

\begin{proof} 
Prove by induction on $n$ that $\mathfrak{M}, w_1 \nmodels I(H)$ iff
there is a branch of worlds $w_1$, \dots, $w_n$ in~$\mathfrak{M}$ with
$w_i R w_{i+1}$ such that for each $i$, $\mathfrak{M}$ is a
countermodel to $\Gamma_i \Seq \Delta_i$ at~$w_i$. This is clear for
$n=1$. Let $H = \Gamma_2 \Seq \Delta_2 \HSeq \dots \HSeq \Gamma_n
\Seq \Delta_n$. Now $\mathfrak{M}, w_1 \nmodels (\bigwedge \Gamma_1
\to \bigvee \Delta_1) \lor \Box I(H)$ iff both $\mathfrak{M}, w_1
\nmodels \bigwedge \Gamma_1 \to \bigvee \Delta_1$ and, for some $w_2$
such that $w_1 R w_2$, $\mathfrak{M}, w_2 \nmodels I(H)$. By induction
hypothesis, the latter holds iff there is a branch of worlds $w_2$,
\dots, $w_n$ such that $\mathfrak{M}$ is a counterexample to $\Gamma_i
\Seq \Delta_i$ at $w_i$ for each $i = 2$, \dots,~$n$. Since
$\mathfrak{M}, w_1 \nmodels \bigwedge \Gamma_1 \to \bigvee \Delta_1$
just means that $\mathfrak{M}$ is a countermodel to $\Gamma_1 \Seq
\Delta_1$ at $w_1$, the claim follows.
\end{proof}

So Parisi's relational hypersequents and Lellman's linear nested
sequents have the same semantic interpretation. Their \emph{calculi}
differ, however, in that Parisi's systems conform to Do\v{s}en's
principle. The base calculus RK contains a pair of rules for~$\Box$,
and extensions of RK for other systems add structural rules, but no
rules that mention~$\Box$. 

Soundness proofs for the relational calculi can be found in
\citet{Parisi2017}. Most of the cases are routine; we give the cases
for $\land$R, $\Box$, EWL and EWR as examples.

\begin{thm}[Soundness]
If $\vdash_\mathrm{RK} H$, then there is no counter-example to $H$.
\end{thm}

\begin{proof}
The proof proceeds by induction on the length of a derivation $\delta$.

\begin{enumerate}
\item Base Case: $\beta$ is an instance of an axiom, ($\varphi \Seq
  \varphi$).  There is no model $\mathfrak{M}$ and possible world
  $w_i$ such that $\mathfrak{M}, w_i \models \varphi$ and
  $\mathfrak{M}, w_i \not \models \varphi$.

\medskip
\quad Let $\beta$ be the last inference of $\delta$. We show that each
rule preserves validity: if the conclusion is not valid, then one of
the premises is not valid. We give the details for $\land$R, $\Box$L,
$\Box$R, and EWR; the other cases are treated the same.
\item $\beta$ is an instance of the $\land$R rule:
\[
\AxiomC{$G \HSeq \Gamma_i \Seq \Delta_i, \phi \HSeq H$}
\AxiomC{$G \HSeq \Gamma_i \Seq \Delta_i, \psi \HSeq H$}
\RightLabel{$\land$R}
\BinaryInfC{$G \HSeq  \Gamma_i \Seq
\Delta_i, \phi \land \psi \HSeq H$}
\DisplayProof
\] 

Let $\mathfrak{M}, w_1,
\ldots, w_n$ form a counter-example to the conclusion. So $w_1,
\ldots, w_n$ forms a branch of worlds such that $w_{k} R w_{k+1}$ and
$\mathfrak{M}$ is a countermodel to each sequent $\Gamma_k \Seq
\Delta_k$ at $w_k$ for all $1 \leq k \leq n$ with $k\neq i$, and a
countermodel to $\Gamma_i \Seq \Delta_i, \phi \land \psi$ at $w_i$. So
$\mathfrak{M}, w_i \models \Gamma_i$, and for each $\theta \in
\Delta_i  \cup \{\phi \land \psi\}$, $\mathfrak{M}, w_i \nmodels
\theta$. It follows that $\mathfrak{M}, w_i \nmodels \phi$ or
$\mathfrak{M}, w_i \nmodels \psi$. In the former case, $\mathfrak{M}$,
$w_1$, \dots,~$w_n$ is a counter-example to the left premise, in the
latter, a counter-example to the right premise.

\item $\beta$ is an instance of $\Box$L. 
\[
\def\fCenter{\HSeq}
\Axiom$G \HSeq  \phantom{\Box \varphi,} \Gamma_{i} \Seq
\Delta_{i} \fCenter \varphi, \Gamma_{i+1} \Seq \Delta_{i+1} \HSeq H$
\RightLabel{$\Box$L}
\UnaryInf$G \HSeq \Box \varphi, \Gamma_{i} \Seq \Delta_{i}
\fCenter \phantom{\varphi,} \Gamma_{i+1} \Seq \Delta_{i+1} \HSeq H$
\DisplayProof
\]

Suppose $\mathfrak{M}, w_1, \ldots, w_n$ is a counter-example for the
conclusion. Then $\mathfrak{M}$ is a countermodel to $\Box \varphi,
\Gamma_i \Seq \Delta_i$ at $w_i$. It follows that $\mathfrak{M}, w_i
\models \Box\varphi$. So, for all  $v$ such that $w_i R v$,
$\mathfrak{M}, v \models \varphi$. Since $w_i R w_{i+1}$, we have that
$\mathfrak{M}, w_{i+1} \models \varphi$. Since $\mathfrak{M}, w_1,
\ldots, w_n$ is a counter-example to the lower hypersequent,
$\mathfrak{M}, w_{i+1} \models \Gamma_{i+1}$ and $\mathfrak{M},
w_{i+1} \nmodels \theta$ for all $\theta \in \Delta_{i+1}$. So
$\mathfrak{M}$ is also a countermodel to the sequent ($\varphi,
\Gamma_{i+1} \Seq \Delta_{i+1}$) at $w_{i+1}$. Hence, $\mathfrak{M},
w_1, \ldots, w_n$ is also a counter-example to the premise.

\item $\beta$ is an instance of $\Box$R.
\[
\def\fCenter{\HSeq}
\Axiom$H \fCenter \Gamma_n \Seq \Delta_n \phantom{\Box \varphi,}
\HSeq \quad \Seq \varphi$
\RightLabel{$\Box$R}
\UnaryInf$H \fCenter \Gamma_n \Seq \Delta_n, \Box \varphi$
\DisplayProof \vspace{3mm}
\]

Suppose that $\mathfrak{M}, w_1, \ldots, w_n$ forms a counter-example
to the conclusion.  Then $\mathfrak{M}$ is a countermodel to $\Gamma_n
\Seq \Delta_n, \Box \varphi$ at~$w_n$.  This means that $\mathfrak{M},
w_n \not\models \Box\varphi$. So, there is some possible world $v$
such that $w_n R v$ and $\mathfrak{M}, v \not\models \varphi$. But
then $\mathfrak{M}$ is also a countermodel to the sequent ($\Seq
\varphi$) at $v$. So $\mathfrak{M}, w_1, \ldots, w_n, v$ is a
counter-example to the premise.

\item $\beta$ is an instance of EWR.
\[\fCenter{{}}
\Axiom$G \fCenter$
\RightLabel{EWR}
\UnaryInf$G \fCenter\HSeq \quad \Seq \quad$
\DisplayProof
\]

Suppose that $\mathfrak{M}, w_1, \ldots, w_n, w_{n+1}$ forms a
counter-example to the conclusion of the inference: If 
\[  
G = \Gamma_1 \Seq \Delta_1 \HSeq \dots \HSeq \Gamma_n \Seq \Delta_n, 
\]$
\mathfrak{M}$ is a countermodel to $\Gamma_k \Seq \Delta_k$ at $w_k$
for all $1 \leq k \leq n$. Of course, $\mathfrak{M}$ is a countermodel
to the empty sequent at any world, in particular~$w_{n+1}$. It follows
that $\mathfrak{M},$ $w_1, \ldots, w_{n}$ is a counter-example to~$G$.
\end{enumerate}
\end{proof}

\begin{thm}
If $\vdash_\mathrm{RT} H$, then there is no reflexive counter-example to $H$.
\end{thm}

\begin{proof}
We show that the EC rule is sound for reflexive frames. Consider:
\[
\AxiomC{$G \HSeq \Gamma_i \Seq \Delta_i \HSeq \Gamma_i
\Seq \Delta_i\HSeq H$}
\RightLabel{EC}
\UnaryInfC{$G \HSeq \Gamma_i \Seq \Delta_i \HSeq H$}
\DisplayProof 
\]

Let $\mathfrak{M}$ be a reflexive counter-example to the conclusion, i.e,
$\mathfrak{M}$ is a countermodel to $\Gamma_k \Seq \Delta_k$ at $w_k$
for $1 \leq k \leq n$. Since the frame is reflexive, $w_i R w_i$. So
$w_1, \ldots, w_i, w_i, \ldots, w_n$ is a branch of worlds where $w_j$
is a counter-example to $\Gamma_j \Seq \Delta_j$ for all $1 \leq j
\leq n$. But this means that $\mathfrak{M}$ is a counter-example to
the premise $G \HSeq \Gamma_i \Seq \Delta_i \HSeq \Gamma_i \Seq
\Delta_i \HSeq H$.
\end{proof}

\begin{thm}
If $\vdash_\mathrm{RB} H$, then there is no symmetric counter-example to $H$.
\end{thm}

\begin{proof}
We show that the Sym rule is sound for symmetric frames. 
\[
\def\fCenter{\ldots}
\Axiom$\Gamma_1 \Seq \Delta_1 \HSeq \fCenter \HSeq
\Gamma_n \Seq \Delta_n$
\RightLabel{Sym} 
\UnaryInf$\Gamma_n \Seq \Delta_n \HSeq \fCenter \HSeq
\Gamma_1 \Seq \Delta_1$
\DisplayProof
\]

Let $\mathfrak{M}, w_1, \ldots, w_n$ be a symmetric counter-example to
the conclusion, i.e., $w_n, \ldots, w_1$ is a branch of worlds such
that $w_{i+1} R w_{i}$ and $\mathfrak{M}$ is a countermodel to
$\Gamma_i \Seq \Delta_i$ at $w_i$ for $1 \leq i \leq n$. Since the
frame is symmetric, it follows that $w_{i} R w_{i+1}$. So $w_1,
\ldots, w_n$ also forms a branch of worlds such that $\mathfrak{M}$ is
a countermodel to each $\Gamma_i \Seq \Delta_i$ at $w_i$ for all $1
\leq i \leq n$, i.e., $\mathfrak{M}$ is a counter-example to the
premise.
\end{proof}

\begin{thm}
If $\vdash_\mathrm{RD} H$, then there is no serial counter-example to $H$.
\end{thm}

\begin{proof}
We show that the Drop rule is sound for serial frames.
\[
\def\fCenter{\ldots}
\Axiom$\Gamma_1 \Seq \Delta_1 \HSeq \fCenter \HSeq
\Gamma_n \Seq \Delta_n \HSeq \quad \Seq\quad$
\RightLabel{Drop} 
\UnaryInf$\Gamma_1 \Seq \Delta_1 \HSeq \fCenter \HSeq
\Gamma_n \Seq \Delta_n$
\DisplayProof
\]
Let $\mathfrak{M}, w_1, \ldots, w_n$ be a serial counter-example to the
conclusion. This means that there is a branch of worlds $w_1, \ldots,
w_n$ such that $w_i R w_{i+1}$ and $\mathfrak{M}$ is a countermodel to
each $\Gamma_i \Seq \Delta_i$ at $w_i$ for all $1 \leq i \leq n$.
Since the frame is serial, there is some world~$v$ such that $w_{n} R
v$. $\mathfrak{M}$ is a countermodel to the empty sequent at~$v$. It
follows that $\mathfrak{M}$, $w_1$, \dots, $w_n$, $v$ is a
counter-example to the premise.
\end{proof}

\begin{thm}
If $\vdash_\mathrm{R4} H$, then there is no transitive counter-example to $H$.
\end{thm}

\begin{proof}
It suffices to show that the EW rule is sound for transitive frames.
\[
\AxiomC{$G \HSeq \Gamma_i \Seq \Delta_i \HSeq \Gamma_{i+1} 
\Seq \Delta_{i+1} \HSeq H$}
\RightLabel{EW}
\UnaryInfC{$G \HSeq \Gamma_i \Seq \Delta_i \HSeq 
\quad \Seq \quad \HSeq \Gamma_{i+1} \Seq \Delta_{i+1} \HSeq
H$}
\DisplayProof
\]

Let $\mathfrak{M}, w_1, \ldots, w_n$ be a transitive counter-example
to the conclusion, where 
\begin{align*}
G & = \Gamma_1 \Seq \Delta_1 \HSeq \ldots \HSeq
      \Gamma_{i-1}\Seq \Delta_{i-1} \text{ and}\\
H & = \Gamma_{i+2} \Seq \Delta_{i+2} \HSeq \ldots \HSeq 
      \Gamma_n \Seq \Delta_n.
\end{align*}
So there is a branch of worlds $w_1, \ldots, w_i, v, w_{i+1}, \ldots,
w_n$ such that $w_k R w_{k+1}$ for all $1 \leq k \leq n$, $w_{i} R v$,
$v R w_{i+1}$ and $\mathfrak{M}$ is a countermodel to each $\Gamma_k
\Seq \Delta_k$ at $w_k$ for all $1 \leq k \leq n$. Since $v R
w_{i+1}$, $w_i R v$, and the frame is transitive, it follows that $w_i
R w_{i+1}$. So $w_1, \ldots, w_i, w_{i+1}, \ldots, w_n$ is also a
branch of worlds and $\mathfrak{M}$ is a countermodel to each
$\Gamma_k \Seq \Delta_k$ at $w_k$ for all $1 \leq k \leq n$. But this
means that $\mathfrak{M}$ is also a counter-example to the premise.

If either $i=0$ (that is, $G \HSeq \Gamma_i \Seq \Delta_i$ is empty)
or $i=n$ ($\Gamma_{i+1} \Seq \Delta_{i+1} \HSeq H$ is empty), then
this is an application of EWL or EWR, respectively, which we have
shown are sound.
\end{proof}

\begin{thm}
If $\vdash_\mathrm{RS5} H$, then there is no reflexive, transitive and symmetric
counter-example to $H$.
\end{thm}

\begin{proof}
We have already shown that EC is sound for all reflexive frames, and
EW for all transitive frames. We now show that the EE rule is sound
for transitive, symmetric frames. Together this means that there
cannot be a reflexive, symmetric, and transitive countermodel to~$H$.

Suppose that $\beta$ is an application of the EE rule.
\[
\def\fCenter{\HSeq}
\AxiomC{$\vdots$}
\noLine
\UnaryInf$G \HSeq \Gamma_i \Seq
\Delta_{i\phantom{+1}} \fCenter \Gamma_{i+1} \Seq \Delta_{i+1}
\HSeq H$
\RightLabel{EE}
\UnaryInf$G \HSeq\Gamma_{i+1} \Seq \Delta_{i+1} \fCenter
\phantom{_+1}\Gamma_i \Seq \Delta_{i\phantom{+1}} \HSeq H$
\DisplayProof
\]
Suppose again that
\begin{align*}
G & = \Gamma_1 \Seq \Delta_1 \HSeq \ldots \HSeq
      \Gamma_{i-1}\Seq \Delta_{i-1} \text{ and}\\
H & = \Gamma_{i+2} \Seq \Delta_{i+2} \HSeq \ldots \HSeq 
      \Gamma_n \Seq \Delta_n
\end{align*}
and that the conclusion has a counterexample, i.e., there is a
symmetric, transitive model $\mathfrak{M}, w_1, \ldots, w_n$ and a
branch of worlds $w_1, \ldots, w_{i-1}, w_{i+1}, w_i, w_{i+2}, \ldots,
w_n$ such that $\mathfrak{M}$ is a countermodel to $\Gamma_k \Seq
\Delta_k$ at $w_k$ for all $1 \leq k \leq n$. We know that
$w_{i-1}Rw_{i+1}$, $w_{i+1} R w_{i}$ and $w_{i} R w_{i+2}$. Since the
frame is transitive, $w_{i-1} R w_{i}$ and $w_{i+1} R w_{i+2}$. Since
the frame is also symmetric, $w_{i} R w_{i+1}$. It follows that $w_1$,
\dots, $w_{i-1}$, $w_i$, $w_{i+1}$, $w_{i+2}$, \dots, $w_n$ is a
branch of worlds in $\mathfrak{M}$. Since $\Gamma_k \Seq \Delta_k$ at
$w_k$ for $1 \leq k \leq n$, there is a counter-example to the
premise. If $G$ or $H$ is empty, the argument is similar.
\end{proof}

\section{Cut-free completeness for RK} \label{sec:comp}

To prove cut-free completeness, we show that for every unprovable
hypersequent~$H$, there is a counter-example. The counter-example is
obtained by constructing a tree $T \subseteq \mathbb{N}^*$ ordered
by a relation~$R$ and assigning labelled hypersequents to
elements~$\sigma \in T$ such that (a) each hypersequent is unprovable,
(b) it is maximal in this respect (``fully reduced''), (c) component
sequents labelled by $\sigma$ of any two hypersequents assigned to
elements of $T$ are identical.  We define a model using $T$, $R$, and
$V$ where $\sigma \in V(p)$ iff $p$ occurs on the left-hand side (lhs)
of any (and thus, by (c), all) component sequents labelled~$\sigma$).
We then show that this model falsifies every component $\Gamma
\Seq[\sigma] \Delta$ at $\sigma$. This relies on the fact that each
hypersequent is fully reduced and on how hypersequents
were assigned to successors of~$\sigma$.

\begin{definition}\label{tree}
  Let $\mathbb{N}^*$ be the set of finite sequences of natural
  numbers. If $\sigma \in \mathbb{N}^*$, then $\sigma.n$ is the
  sequence~$\sigma$ extended by~$n \in \mathbb{N}$.

  A subset~$T$ of $\mathbb{N}^*$ is a \emph{tree} iff whenever
  $\sigma.n \in T$ then $\sigma \in T$. We consider four relations
  on~$\mathbb{N}^*$ and~$T$: 
  \begin{enumerate}
  \item The successor relation~$R^1$: $\sigma R^1 \sigma'$ iff
    $\sigma' = \sigma.n$.
  \item The reflexive closure~$R^=$ of $R$.
  \item The transitive closure~$R^+$ of $R$.
  \item The reflexive transitive closure $R^*$ of~$R$. 
  \end{enumerate}
\end{definition}

Obviously $R^=$ is reflexive, $R^+$ is transitive, and $R^*$ is
reflexive and transitive, both on $\mathbb{N}^*$ and any tree~$T$.

\begin{definition}\label{branch}
  A sequence $\Sigma = \sigma_1$, \dots, $\sigma_n$ is an
  $R$-\emph{branch} iff $\sigma_i R \sigma_{i+1}$ for $1 \le i
  < n$. $\Sigma$ is an $R$-\emph{path} through $T$ if it is an
  $R$-branch of $T$ and for all $\sigma \in T$, not $\sigma R^1
  \sigma_1$ and not $\sigma_n R^1 \sigma$ (i.e., $R$-paths are
  $R$-branches that are maximal with respect to~$R^1$, although they
  need not be maximal in the order~$R$).
\end{definition}

\begin{definition}
  An \emph{$R$-labelled hypersequent} is a sequence $\Gamma_1
  \Seq[\sigma_1] \Delta_1 \HSeq \dots \HSeq \Gamma_n \Seq[\sigma_n]
  \Delta_n$ where $\sigma_1$, \dots, $\sigma_n$ is an $R$-branch.

If $H$ is an $R$-labelled hypersequent, then $H(\sigma)$ is the rightmost  component sequent $\Gamma \Seq[\sigma] \Delta$ of~$H$ or the empty
  sequent if $H$ has no such component sequent.
  
  $\Sigma(H)$ is the sequence of labels $\sigma_1,\dots,\sigma_n$ of
  the component sequents of~$H$.
  
  We write $\Gamma(H,\sigma)$ and $\Delta(H,\sigma)$ for the left-hand
  side and right-hand side of $H(\sigma)$, i.e., the sets of formulas
  such that $H(\sigma)$ is $\Gamma(H, \sigma) \Seq[\sigma] \Delta(H,
  \sigma)$
  
  We say a sequent $\Gamma' \Seq \Delta'$ \emph{extends} a sequent $\Gamma
  \Seq \Delta$ iff $\Gamma \subseteq \Gamma'$ and $\Delta \subseteq
  \Delta'$. If $H$ and $H'$ are labelled hypersequents, we say $H'$
  \emph{extends} $H$ iff for all $\sigma$ occurring as labels in $H$,
  $H'(\sigma)$ extends $H(\sigma)$.
\end{definition}

Our construction will produce a tree of labels~$\sigma$ and unprovable
$R^1$-labelled hypersequents. From this tree we will extract a
counter-example. We will ensure that new hypersequents added to the
tree are always extensions of original ones. In the next section, we
will extend the construction to RT, in which case we deal with
$R^+$-labelled hypersequents which may contain more than one component
sequent with the same label~$\sigma$. However, the construction will
guarantee that if a hypersequent contains two component sequents
$\Gamma \Seq[\sigma] \Delta$ and $\Gamma' \Seq[\sigma] \Delta$ with
the same label~$\sigma$, the component further to the right extends
the component to the left in~$H$, by defining reducts (in the
following definition) always on the basis of the rightmost component
labelled by~$\sigma$.  The definition of the model, specifically, the
valuation at~$\sigma$, then also need only take into account the
rightmost component~$H(\sigma)$. For the remainder of this section,
however, we will deal with $R^1$-labelled hypersequents only.

\begin{definition}\label{reduct}
  Given an $R$-labelled hypersequent~$H$ and a label~$\sigma$, we define a
  \emph{$\sigma$-reduct} (corresponding to a rule) of $H$ as the
  corresponding hypersequent on the right in Table~\ref{reducts}.
\begin{table}[h]
\caption{Reducts of labelled hypersequents}
\label{reducts}
\centering
  \begin{tabular}{l@{\quad}l@{\quad}l}
  \hline\hline
    Rule & Hypersequent $H$ & $\sigma$-Reduct of $H$\\ \hline
    $\lnot$L & $G \HSeq \lnot \phi, \Gamma \Seq[\sigma] \Delta \HSeq G'$ &
    $G \HSeq \lnot \phi, \Gamma \Seq[\sigma] \Delta, \phi \HSeq G'$ \\
    $\lnot$R & $G \HSeq \Gamma \Seq[\sigma] \Delta, \lnot \phi \HSeq G'$ &
    $G \HSeq \phi, \Gamma \Seq[\sigma] \Delta, \lnot\phi \HSeq G'$ \\
$\land$L & $G \HSeq \phi \land \psi, \Gamma \Seq[\sigma] \Delta \HSeq G'$ &$G \HSeq \phi, \psi, \phi \land \psi, \Gamma \Seq[\sigma] \Delta \HSeq G'$
\\
$\land$R & $G \HSeq \Gamma \Seq[\sigma] \Delta, \phi \land \psi \HSeq G'$ &    $G \HSeq \Gamma \Seq[\sigma] \Delta, \phi \land \psi, \phi\HSeq G'$ \\
    & & if unprovable, otherwise\\
& & $G \HSeq \Gamma \Seq[\sigma] \Delta, \phi \land \psi, \psi\HSeq G'$ \\$\Box$L & $G \HSeq \Box \phi, \Gamma' \Seq[\sigma'] \Delta' \HSeq \Gamma
\Seq[\sigma] \Delta \HSeq G'$ &
    $G \HSeq \Box\phi, \Gamma' \Seq[\sigma'] \Delta' \HSeq \phi,
    \Gamma \Seq[\sigma] \Delta \HSeq G'$ \\
      \hline\hline
  \end{tabular}
  \end{table}

  In each case, the displayed component sequent labelled by $\sigma$
  is the rightmost such in~$H$, if there is more than one. 
  
  A hypersequent is called \emph{$\sigma$-reduced} if
  it is identical to all of its $\sigma$-reducts, otherwise it is
  \emph{$\sigma$-reducible}.  If it is $\sigma$-reduced for all $\sigma$
  occurring in it as labels, it is called \emph{fully reduced}.
\end{definition}

\begin{prop}\label{reduct-unprovable}
  If $H$ is unprovable, any $\sigma$-reduct of it is also unprovable.
\end{prop}

\begin{proof}
  If the $\sigma$-reduct of $H$ were provable, the relevant rule
  would prove~$H$. For instance, suppose
  $G \HSeq \Gamma \Seq[\sigma] \Delta, \phi \land \psi \HSeq G'$ is
  unprovable. Then one of
  \begin{align*}
    & G \HSeq \Gamma \Seq[\sigma] \Delta, \phi \land \psi, \phi\HSeq G'\\
    & G \HSeq \Gamma \Seq[\sigma] \Delta, \phi \land \psi, \psi\HSeq G'
  \end{align*}
  must be unprovable. For suppose both were provable. Then we'd have:
  \begin{prooftree}
    \def\fCenter{\HSeq}
  \AxiomC{}
\DeduceC{$G \HSeq \Gamma \Seq[\sigma] \Delta, \phi \land \psi, \phi\HSeq
G'$}
  \AxiomC{}
\DeduceC{$G \HSeq \Gamma \Seq[\sigma] \Delta, \phi \land \psi, \psi\HSeq
G'$}
  \RightLabel{$\land$R}
\BinaryInf$G \fCenter \Gamma \Seq[\sigma] \Delta, \phi \land \psi \HSeq G'$\end{prooftree}
  Or, suppose the $\sigma$-reduct based on the $\Box$L-rule were
  provable. Then we'd have:
  \[\bottomAlignProof
    \def\fCenter{\HSeq}
    \AxiomC{}
\Deduce$G \HSeq \Box\phi, \Gamma' \Seq[\sigma'] \Delta \fCenter \phi,
\Gamma \Seq[\sigma] \Delta \HSeq G'$
    \RightLabel{$\Box$L}
    \UnaryInf$G \HSeq \Box\phi, \Gamma' \Seq[\sigma'] \Delta \fCenter
    \Gamma \Seq[\sigma] \Delta \HSeq G'$ 
  \DisplayProof\]
\end{proof}

\begin{prop}\label{full-red-unprov}
  Every unprovable labelled hypersequent~$H$ is extended by an
  unprovable, fully reduced hypersequent~$\Red(H)$ (called its
  \emph{full reduction}).
\end{prop}

\begin{proof}
  If $H$ is already fully reduced, we have nothing to prove.
  Otherwise, there is a least $\sigma$ (in the 
  prefix order~$R^*$) so that $H$ is not $\sigma$-reduced. Any
  $\sigma$-reduction of a reducible hypersequent extends it. So,
  starting with the set $\{H\}$ and adding $\sigma$-reductions results
  in a set of unprovable hypersequents, partially ordered by
  extension. This set is finite, as can easily be seen by induction on
  the number and degree of formulas in $H(\sigma)$ and the number of
  formulas of the form $\Box \phi$ in $H(\sigma')$.  A maximal element
  in this order is an unprovable $\sigma$-reduced hypersequent
  extending~$H$. The proposition follows by induction on the number of
  components of~$H$.
\end{proof}

\begin{prop}\label{red-subf}
  Let $H$ be $\Red(H')$ for some unprovable sequent~$H'$, and let
  $H(\sigma)$ be $\Gamma \Seq[\sigma] \Delta$.
  \begin{enumerate}
  \item\label{rs-nl} If $\lnot \phi \in \Gamma$, then $\phi
    \in \Delta$.
  \item\label{rs-nr} If $\lnot \phi \in \Delta$, then $\phi
    \in \Gamma$.
  \item\label{rs-al} If $\phi \land \psi \in \Gamma$, then
    $\phi \in \Gamma$ and $\psi \in \Gamma$.
  \item\label{rs-ar} If $\phi \land \psi \in \Delta$, then
    $\phi \in \Delta$ or $\psi \in \Delta$.
  \item\label{rs-box} If $\Box \phi \in \Gamma$, $\sigma R^1 \tau$, and
    $\tau$ occurs in $H$, then $\phi \in \Gamma(H,\tau)$.
  \end{enumerate}
\end{prop}

\begin{proof}
  Since $H$ is $\sigma$-reduced, $H(\sigma)$ is identical to all its
$\sigma$-reducts. Inspection of the definition of reducts
(Table~\ref{reducts}) establishes (\ref{rs-nl})--(\ref{rs-ar}).

  For (\ref{rs-box}), suppose that $\Box \phi \in \Gamma$,
  $\sigma R^1 \tau$ and $\tau$ occurs as a label in~$H$. Since $H$
  is $\sigma$-reduced, $H$ is identical to its $\Box$L
  $\sigma$-reducts. Since $\Sigma(H)$ is an $R^1$-branch, the component
  $H(\tau)$ occurs immediately to the right of $H(\sigma)$, i.e.,
  $\Gamma(H,\tau) \ni \phi$.
\end{proof}

\begin{definition}\label{successor}
  Suppose $H=G \HSeq \Gamma \Seq[\sigma] \Delta \HSeq G'$ is an
  unprovable fully reduced hypersequent, and $\Box \psi \in \Delta$. The
\emph{$\sigma.n$-$\psi$-successor}~$\Succ_{\sigma.n}^\psi(H)$ of~$H$ is the  hypersequent $\Red(G \HSeq \Gamma \Seq[\sigma] \Delta \HSeq \quad
  \Seq[\sigma.n] \psi)$.
\end{definition}

We record some facts about the successor construction.

\begin{prop}\label{succ-unprov}
\ 
  \begin{enumerate}
  \item\label{su-unprov} The $\sigma.n$-$\psi$-successor of an
    unprovable fully reduced hypersequent is unprovable.
  \item\label{su-unchanged} If $H^*$ is a $\sigma.n$-$\psi$-successor
    of~$H$ and $\tau R^* \sigma$, then $H(\tau) = H^*(\tau)$
    (i.e., passing to successors does not change the sequent labelled
    $\sigma$ or any to the left of it).
  \end{enumerate}
\end{prop}

\begin{proof}
Suppose $H$ is a fully reduced hypersequent of the form 
  \begin{align*}
  H& =G \HSeq \Gamma
  \Seq[\sigma] \Delta, \Box \psi \HSeq G' \text{ and}\\
  H'& =G \HSeq \Gamma
  \Seq[\sigma] \Delta, \Box \psi \HSeq \quad
  \Seq[\sigma.n] \psi.
  \end{align*}
  \begin{enumerate}
  \item If $H'$ were provable, then $H$ would be provable:
    \begin{prooftree}
    \def\fCenter{\Seq}
    \AxiomC{}
\Deduce$G \HSeq \Gamma \fCenter \Delta, \Box \psi \HSeq \quad \Seq \psi $    \RightLabel{$\Box$ R}
    \UnaryInf$G \HSeq \Gamma \fCenter \Delta, \Box \psi$
    \RightLabel{EWR}
    \UnaryInf$G \HSeq \Gamma \fCenter \Delta, \Box\psi \HSeq G'$
  \end{prooftree}
  $\Succ_{\sigma.n}^\psi(H)$ is $\Red(H')$, which is unprovable if~$H'$
  is by Proposition~\ref{full-red-unprov}.

  \item Since $H$ is fully reduced, $\Seq[\sigma.n] \psi$ is the only
    reducible sequent in $H'$. So, for every $\tau$ in $\Sigma(H')$
    other than~$\sigma.n$, every $\tau$-reduct of $H'$ is identical to
    $H'$. Moreover, in constructing $\tau$-reducts, no formulas are
    added to component sequents to the left of $H'(\tau)$. In
    particular, reduction of $H'(\sigma.n)$ does not affect
    $H'(\sigma)$, throughout the construction of~$\Red(H')$ given in
    the proof of Proposition~\ref{full-red-unprov}. Hence, for all
    $\tau$ with $\tau R^* \sigma$, $H(\tau) = \Red(H')(\tau) =
    H^*(\tau)$.
  \end{enumerate}
\end{proof}

\begin{definition}\label{succ-tree}
  Let $H$ be an unprovable hypersequent
  \[
  \Gamma_1 \Seq \Delta_1 \HSeq \dots \HSeq \Gamma_n \Seq \Delta_n.
  \]
  Let $H'$ be the full reduction of
  \[
\Gamma_1 \Seq[\sigma_1] \Delta_1 \HSeq \dots \HSeq \Gamma_n \Seq[\sigma_n]
\Delta_n
  \]
  with $\sigma_i = 0\dots0$ with $i$ $0$'s, and let $H'_i = \Gamma_1
  \Seq[\sigma_1] \Delta_1 \HSeq \dots \HSeq \Gamma_i \Seq[\sigma_i]
  \Delta_i$ ($i \le n$).

  We define a partial mapping~$\lambda$ from $\mathbb{N}^*$ to
  labelled hypersequents inductively. Assuming $\lambda(\sigma)$ is
  already defined, let $\Gamma \Seq[\sigma] \Delta$ be
  $\lambda(\sigma)(\sigma)$, i.e., the (rightmost) $\sigma$-labelled
  component of the hypersequent $\lambda(\sigma)$, and let $\psi_1$,
  \dots, $\psi_l$ be all the formulas $\psi_k$ such that $\Box \psi_k \in
  \Delta$.
  \begin{align*}
    \lambda(0) & = H_1' \\
    \lambda(\sigma.k) & = \begin{cases}
      H_{i+1}' & \text{if $k=0$, $i < n$, and $\sigma=\sigma_{i}$}\\
      \Succ_{\sigma.k}^{\psi_k}(\lambda(\sigma)) &
      \text{if $k>0$, $\lambda(\sigma)$ is defined, and $\psi_k$ exists}\\
      \text{undefined} & \text{otherwise}
      \end{cases}
  \end{align*}

  Let $T$ be the set of all $\sigma \in \mathbb{N}^*$ such that
  $\lambda(\sigma)$ is defined and let $S(H) = \{\lambda(\sigma) : \sigma
  \in T\}$ be all labelled hypersequents in the range of~$\lambda$.
\end{definition}

\begin{prop}\label{sigma-tree}
  We record some facts about $T$, $\lambda$, and $S(H)$:
  \begin{enumerate}
    \item\label{st-tree} $T$ is a tree.
\item\label{st-unprov} If $G \in S(H)$, $G$ is unprovable and fully
reduced.
\item\label{st-sigma} If $G = \lambda(\sigma)$, $\tau$ occurs in $G$ iff      $\tau R^*\sigma$.
    \item\label{st-unchanged} If $G = \lambda(\sigma)$, $G' =
\lambda(\tau)$, and $\sigma R^* \tau$, then $G(\sigma) = G'(\sigma)$.    \item\label{st-equal} If $G$, $G' \in S(H)$ and $\sigma$ occurs in
      both, $G(\sigma) = G'(\sigma)$.
    \item\label{st-successor} If $H' = G \HSeq \Gamma \Seq[\sigma]
      \Delta \HSeq G' \in S(H)$ and $\Box \psi \in \Delta$, there is a
      $k$ such that $\sigma.k \in T$ such that $H'' = G \HSeq
      \Gamma \Seq[\sigma] \Delta \HSeq \Gamma' \Seq[\tau] \Delta'
      \in H(S)$ and $\psi \in \Delta'$. 
  \end{enumerate}
\end{prop}

\begin{proof}
  \begin{enumerate}
    \item By construction, if $\lambda(\sigma)$ is undefined,
      $\lambda(\sigma.n)$ is undefined. Hence, if $\sigma \in T$ has
      the property that $\tau \in T$ for all $\tau R^+ \sigma$, so
      does $\sigma.n \in T$.
  
    \item Each $\lambda(\sigma_i)$, i.e., $H'_i$ for $i = 1$,
      \dots,~$n$, is fully reduced and unprovable (If $H'_i$ is
      provable, so is $H' = H'_i \HSeq G$, by EWR.)  By induction on
      $\sigma$, and Propositions \ref{full-red-unprov} and
      \ref{succ-unprov}(\ref{su-unprov}), each $\lambda(\sigma.k)$
      ($\sigma.k \in T$ and $k>0$) is unprovable (and fully reduced by
      construction).

    \item By induction on $\sigma \in T$. There is no $\tau$ such that
      $\tau R^1 0$. The property holds for $\sigma_i$ by definition of
      $\lambda(\sigma_i)$. It holds for $\sigma.k$ ($k \ge 1$) by the
      definition of $\lambda(\sigma.k)$ and
      Proposition~\ref{succ-unprov}(\ref{su-unchanged}).
  
    \item By induction on $\sigma \in T$: The definition of $H'$
      ensures the property holds for $\lambda(\sigma_i)$, and the
      definition of $\Succ_{\sigma.k}^\psi$ ensures that if it holds
      for $\lambda(\sigma)$ it also does for $\lambda(\sigma.k)$ ($k
      \ge 1$).

    \item Let $\tau$, $\tau'$ be such that $\lambda(\tau) = G$ and
      $\lambda(\tau') = G'$. If $\sigma$ occurs in both $G$ and $G'$,
      by (\ref{st-sigma}), $\sigma R^* \tau$ and $\sigma R^* \tau'$.
      Let $G'' = \lambda(\sigma)$.  Then by (\ref{st-unchanged}),
      $G''(\sigma) = G(\sigma)$ and $G''(\sigma) = G'(\sigma)$, and so
      $G(\sigma) = G'(\sigma)$.

\item $H''$ is a $\sigma.k$-$\psi$-successor of $H'$.
\end{enumerate}
\end{proof}

Since $G(\sigma) = G'(\sigma) = \Gamma \Seq[\sigma] \Delta$ for any
two $G$, $G' \in S(H)$ which both contain $\sigma$, we can define
$\Gamma(\sigma) = \Gamma$ and $\Delta(\sigma) = \Delta$ independently
of the individual hypersequents in $S(H)$. 

Given an unprovable hypersequent $H$, let $T$ and $S(H)$ be as in
Definition~\ref{succ-tree} and let $\mathfrak{M} = \langle T, R^1,
V\rangle$ with $\sigma \in V(p)$ iff $p \in \Gamma(\sigma)$.

\begin{prop}\label{truth-lemma}
  For all $\phi$, if $\phi \in \Gamma(\sigma)$ then $\mathfrak{M},
  \sigma \models \phi$ and if $\phi \in \Delta(\sigma)$, then
  $\mathfrak{M}, \sigma \nmodels \phi$.
\end{prop}

\begin{proof}
  By induction on $\phi$.

  If $p \in \Gamma(\sigma)$, then $\sigma \in V(p)$ by definition, so
  $\mathfrak{M}, \sigma \models p$.

  If $p \in \Delta(\sigma)$, then $p \notin \Gamma(\sigma)$ (otherwise
  $\Gamma(\sigma) \Seq \Delta(\sigma)$ and any hypersequent containing
  it would be provable.) So $\sigma \notin V(p)$.

  If $\lnot \phi \in \Gamma(\sigma)$, by
  Proposition~\ref{red-subf}(\ref{rs-nl}), $\phi \in \Delta(\sigma)$.
  By induction hypothesis, $\mathfrak{M}, \sigma \nmodels \phi$, so
  $\mathfrak{M}, \sigma \models \lnot \phi$. Similarly for $\lnot \phi
  \in \Delta(\sigma)$, using Proposition~\ref{red-subf}(\ref{rs-nr}).

  If $\phi \land \psi \in \Gamma(\sigma)$, by
  Proposition~\ref{red-subf}(\ref{rs-al}), $\phi \in \Gamma(\sigma)$
  and $\psi \in \Gamma(\sigma)$. By induction hypothesis,
  $\mathfrak{M}, \sigma \models \phi$ and $\mathfrak{M}, \sigma
  \models \psi$, so $\mathfrak{M}, \sigma \models \phi \land \psi$.

  If $\phi \land \psi \in \Delta(\sigma)$, by
  Proposition~\ref{red-subf}(\ref{rs-ar}), $\phi \in \Delta(\sigma)$
  or $\psi \in \Delta(\sigma)$. By induction hypothesis,
  $\mathfrak{M}, \sigma \nmodels \phi$ or $\mathfrak{M}, \sigma
  \nmodels \psi$, so $\mathfrak{M}, \sigma \nmodels \phi \land \psi$.

  Suppose $\Box \phi \in \Gamma(\sigma)$ and let $\sigma R^1 \tau$. By
  Proposition~\ref{sigma-tree}(\ref{st-sigma}) and
  Proposition~\ref{red-subf}(\ref{rs-box}), $\phi \in
  \Gamma(\tau)$. By induction hypothesis, $\mathfrak{M}, \tau \models
  \phi$. Thus, $\mathfrak{M}, \sigma \models \Box \phi$.

  Suppose $\Box \phi \in \Delta(\sigma)$. By
  Proposition~\ref{sigma-tree}(\ref{st-successor}), there is a $\tau$
  such that $\sigma R^1 \tau$ (namely, $\tau = \sigma.k$ for some $k$)
and $\phi \in \Delta(\tau)$. By induction hypothesis, $\mathfrak{M}, \tau  \nmodels \phi$, hence $\mathfrak{M}, \sigma \nmodels \Box \phi$.
\end{proof}

\begin{cor}
  The calculus RK is complete for $K$.
\end{cor}

\begin{example}\label{example}
Consider the hypersequent $\Box\lnot(p \land q) \Seq \Box \lnot q \HSeq p
\Seq$. The counter-example construction begins by labelling the components
using the branch $0$, $0.0$:
  \begin{align*}
  \Box\lnot(p \land q) \Seq[0] \Box \lnot q & \HSeq p \Seq[0.0]
\intertext{It is $0$-reduced, but not $0.0$-reduced. A $0.0$-reduct, using
the
$\Box$L rule, is:}
\Box\lnot(p \land q) \Seq[0] \Box \lnot q & \HSeq \lnot(p \land q), p
\Seq[0.0]
\intertext{In turn, we can apply a $\lnot$L-reduction to the sequent
labelled $0.0$ to obtain}
\Box\lnot(p \land q) \Seq[0] \Box \lnot q & \HSeq \lnot(p \land q), p
\Seq[0.0] p \land q
\intertext{Finally, we apply a $\land$R-reduction to obtain}
\Box\lnot(p \land q) \Seq[0] \Box \lnot q & \HSeq \lnot(p \land q), p
\Seq[0.0] p \land q, q
\intertext{Since $\Box \lnot q \in \Delta(0)$, there is a $0.1$-$\lnot
q$-successor, namely}
\Box\lnot(p \land q) \Seq[0] \Box \lnot q & \HSeq \quad \Seq[0.1] \lnot q\intertext{Its full reduction is}
\Box\lnot(p \land q) \Seq[0] \Box \lnot q & \HSeq q, \lnot(p \land q)
\Seq[0.1] \lnot q, p\land q, p
  \end{align*}
  We now have $T = \{0, 0.0, 0.1\}$ with
  \begin{align*}
    \lambda(0) & =   \Box\lnot(p \land q) \Seq[0] \Box \lnot q \\
\lambda(0.0) & = \Box\lnot(p \land q) \Seq[0] \Box \lnot q \HSeq \lnot(p
\land q), p \Seq[0.0] p \land q, q\\
\lambda(0.1) & = \Box\lnot(p \land q) \Seq[0] \Box \lnot q \HSeq q, \lnot(p
\land q) \Seq[0.1] \lnot q, p\land q, p
  \end{align*}
  The corresponding counter-example is
  
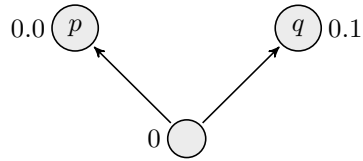
\begin{center}
    \begin{tikzpicture}[modal]
      \node[world] (w1) [label={left:$0$}] {}; 
      \node[world] (w2) [label={left:$0.0$},
        above left=of w1] {$p$}; 
      \node[world] (w3) [label={right:$0.1$},
        above right=of w1] {$q$};
      \draw[->] (w1) to (w2);
      \draw[->] (w1) to (w3);
    \end{tikzpicture}
  \end{center}
\end{example}

\section{Cut-free completeness of RT and RD} \label{sec:comp2}

The completeness proof above can be extended to RT. First we extend
Definition~\ref{reduct} to include the following:
\begin{center}
  \begin{tabular}{l@{\quad}l@{\quad}l}
    Rule & Hypersequent~$H$ & $\sigma$-Reduct of~$H$\\
    EC & $G \HSeq  \Box \phi, \Gamma \Seq[\sigma]  \Delta \HSeq G'$ &
$G \HSeq \Box\phi, \Gamma \Seq[\sigma] \Delta \HSeq \phi, \Box \phi, \Gamma
\Seq[\sigma] \Delta \HSeq G'$
  \end{tabular}
\end{center}
where the sequent $\Box \phi, \Gamma \Seq[\sigma] \Delta$ is
the rightmost sequent labelled~$\sigma$ in~$H$ and $\phi \notin \Gamma$.

Then Proposition~\ref{reduct-unprovable} still holds, i.e., if $H$ is
unprovable, so are its $\sigma$-reducts. 

\begin{prooftree}
  \AxiomC{$G \HSeq \Box\phi, \Gamma \Seq[\sigma] \Delta \HSeq
    \phi, \Box\phi, \Gamma \Seq[\sigma] \Delta \HSeq G'$}
  \RightLabel{$\Box$L}
  \UnaryInfC{$G \HSeq \Box\phi, \Gamma \Seq[\sigma] \Delta \HSeq
    \Box\phi, \Gamma \Seq[\sigma] \Delta \HSeq G'$}
  \RightLabel{EC}
  \UnaryInfC{$G \HSeq \Box\phi, \Gamma \Seq[\sigma] \Delta \HSeq G'$}
\end{prooftree}

Proposition~\ref{full-red-unprov} also still holds for the extended
definition of ``fully reduced'' now including reducts for rule~EC. A
$\sigma$-reduct~$H'$ of a hypersequent $H$ also extends it: Suppose
$\Box \phi, \Gamma$ is $\Gamma(H, \sigma)$, the left side of the
rightmost sequent in~$H$ labelled~$\sigma$. Then $\Gamma(H',\sigma)$,
the left side of the rightmost sequent labelled~$\sigma$ in~$H'$, is
$\phi, \Box \phi, \Gamma$. Clearly, the number of times an EC reduction
can be applied to the sequent labelled~$\sigma$ is bounded by the sum
of the degrees of the formulas in~$H$.

Proposition~\ref{red-subf}(\ref{rs-box}) now holds in the form:
If $H$ is $\Red(H')$ for some hypersequent~$H'$, and $H(\sigma) =
\Gamma \Seq \Delta$, then
\begin{enumerate}
\item[\ref{rs-box}$'$] If $\Box \phi \in \Gamma$, $\sigma
    R^=\tau$, and $\tau$ occurs in $H$, then $\phi \in
    \Gamma(H,\tau)$.
\end{enumerate}
If $\sigma R^1 \tau$, then we just have a case of
Proposition~\ref{red-subf}(\ref{rs-box}). For the case $\sigma=\tau$,
we have to show that if $\Box \phi \in \Gamma$, then $\phi \in
\Gamma$. This holds since $H$ is fully reduced, and $G \HSeq \phi,
\Gamma \Seq[\sigma] \Delta \HSeq G'$ is a $\sigma$-reduct of $G \HSeq
\Gamma \Seq[\sigma] \Delta \HSeq G'$ (for rule EC).

Definition~\ref{succ-tree} yields a tree of unprovable hypersequents
$S(H)$ for any unprovable hypersequent~$H$ also when EC-reductions are
included in the definition of $\Red$. For the definition of the
$\lambda(\sigma.k)$, note that $\lambda(\sigma)(\sigma) = \Gamma
\Seq[\sigma] \Delta$ is the \emph{rightmost} $\sigma$-labelled
component of the hypersequent $\lambda(\sigma)$. Thus, successors are
computed from the fully reduced hypersequent component.

Proposition~\ref{sigma-tree} still holds since it is independent of
the definition of reduction.

Completeness for reflexive models now follows: If $H$ is unprovable,
$S(H)$ is a tree of fully reduced unprovable hypersequents. Define
$\mathfrak{M} = \langle T, R^=, V\rangle$ as before, with the
difference that the accessibility relation is the reflexive
closure~$R^=$ of~$R$. Proposition~\ref{truth-lemma} holds for $S(H)$
and $\mathfrak{M}$, since the only relevant difference is the case
$\Box \phi \in \Gamma(\sigma)$, which holds by
Proposition~\ref{red-subf}(\ref{rs-box}$'$).

\begin{example}\label{ex2}
  Consider the hypersequent $\Box\lnot(p \land q), p \Seq \Box
  \lnot q \HSeq p \Seq$.   
  Again we begin by labelling the components using the branch $0$, $0.0$:
  \begin{align*}
  \Box\lnot(p \land q), p \Seq[0] \Box \lnot q & \HSeq p \Seq[0.0]
\intertext{This hypersequent is not $0$-reduced. A $0$-reduct using the RT
reduction is:}
\Box\lnot(p \land q), p \Seq[0] \Box \lnot q & \HSeq \lnot(p \land q),
\Box\lnot(p \land q), p \Seq[0] \Box \lnot q \HSeq p \Seq[0.0]
  \intertext{which further reduces to}
\Box\lnot(p \land q), p \Seq[0] \Box \lnot q & \HSeq \lnot(p \land q),
\Box\lnot(p \land q), p \Seq[0] \Box \lnot q, p \land q, q \HSeq p
\Seq[0.0]
   \intertext{This is now $0$-reduced. The full reduct, as before, is:}
\Box\lnot(p \land q) \Seq[0] \Box \lnot q & \HSeq \lnot(p \land q),
\Box\lnot(p \land q), p \Seq[0] \Box \lnot q, p \land q, p \HSeq \lnot(p
\land q), p \Seq[0.0] p \land q, q
\intertext{There again is a $0.1$-$\lnot q$-successor, namely}
\Box\lnot(p \land q) \Seq[0] \Box \lnot q & \HSeq \lnot(p \land q),
\Box\lnot(p \land q), p \Seq[0] \Box \lnot q, p \land q, q \HSeq \quad
\Seq[0.1] \lnot q
\intertext{Its full reduction is}
\Box\lnot(p \land q) \Seq[0] \Box \lnot q & \HSeq q, \lnot(p \land q),
\Box\lnot(p \land q) \Seq[0] \Box \lnot q, p \land q, p \HSeq q, \lnot(p
\land q) \Seq[0.1] \lnot q, p\land q, p
  \end{align*}
  The corresponding counter-example is
  
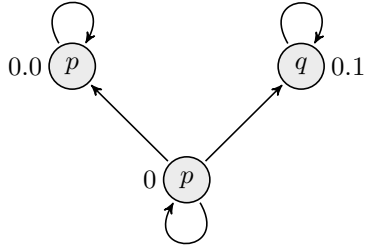
\begin{center}
    \begin{tikzpicture}[modal]
      \node[world] (w1) [label={left:$0$}] {$p$}; 
      \node[world] (w2) [label={left:$0.0$},
        above left=of w1] {$p$}; 
      \node[world] (w3) [label={right:$0.1$},
        above right=of w1] {$q$};
      \draw[->] (w1) to (w2);
      \draw[->] (w1) to (w3);
      \draw[reflexive below] (w1) to (w1);
      \draw[reflexive above] (w2) to (w2);
      \draw[reflexive above] (w3) to (w3);
    \end{tikzpicture}
  \end{center}
\end{example}

To prove completeness of RD for serial models, we have to ensure that
the accessibility relation on~$T$ is serial. To do this, we extend
Definition~\ref{successor}: Suppose $H=G \HSeq \Gamma \Seq[\sigma]
\Delta$ is an unprovable fully reduced hypersequent, and $\Delta$
contains no formula of the form $\Box \psi$ (i.e., it has no
$\sigma.n$-$\psi$ successor, where $\sigma$ is the label of the
rightmost sequent in~$H$).  The
\emph{$\sigma.n$-successor}~$\Succ_{\sigma.n}(H)$ of $H$ is the
hypersequent $\Red(G \HSeq \Gamma \Seq[\sigma] \Delta \HSeq \quad
\Seq[\sigma.n] \quad)$.

Proposition~\ref{succ-unprov} also holds for $\sigma.n$-successors, by
the Drop rule. Definition~\ref{succ-tree} is extended by including the
$\sigma.n$-successor of $H$ if there is no
$\sigma.n$-$\psi$-successor. Proposition~\ref{red-subf}(\ref{rs-box})
still holds since the $\sigma.n$-successor is fully reduced as in the
case for~K. Proposition~\ref{sigma-tree} and~\ref{truth-lemma} still
hold.  The relation $R$ on $T$ in this case is serial, since for every
$\sigma$ that occurs as a label on a sequent either $\sigma.0$ also
occurs as a label in the original labelled hypersequent~$H'$, or there
is a $\sigma.n$-$\psi$ successor, or $\sigma$ is the label of the
rightmost sequent without a formula of the form $\Box \psi$ in the
succedent, in which case there is a $\sigma.n$-successor.

This method of adding successors results in an infinite tree, but we
can do a bit better: only add a $\sigma.n$-successor if $\Box\phi \in
\Gamma(\sigma)$, and add $\langle \sigma,\sigma\rangle$ to the
accessibility relation. For instance, suppose we start with $\Box\Box
p \Seq[0] \Box p$. This is completely reduced, and has a $0.1$-$p$
successor which reduces to
\begin{align*}
  \Box\Box p \Seq[0] \Box p & \HSeq \Box p \Seq[0.1] p
\intertext{There is no $\Box \psi \in \Delta(0.1)$, so a $0.1.1$-successor
is}
\Box\Box p \Seq[0] \Box p & \HSeq \Box p \Seq[0.1] p \HSeq \quad
\Seq[0.1.1]
  \intertext{which reduces to}
  \Box\Box p \Seq[0] \Box p & \HSeq \Box p \Seq[0.1] p \HSeq p \Seq[0.1.1]
\end{align*}
The serial counter-example is
  
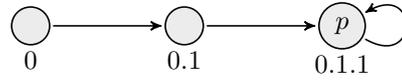
\begin{center}
    \begin{tikzpicture}[modal]
      \node[world] (w1) [label={below:$0$}] {}; 
      \node[world] (w2) [label={below:$0.1$},
        right=of w1] {}; 
      \node[world] (w3) [label={below:$0.1.1$},
        right=of w2] {$p$};
      \draw[->] (w1) to (w2);
      \draw[->] (w2) to (w3);
      \draw[reflexive right] (w3) to (w3);
    \end{tikzpicture}
  \end{center}

\section{Failure of the method for B and 4} \label{sec:fail}

The constructions of the K, T, and D counter-examples to RK, RT, and
RD-unprovable hypersequents work because once a sequent with label
$\sigma$ is reduced, it remains unchanged in the reduction of
successors. This guarantees that in the entire tree of hypersequents,
all (rightmost, in the case of RT) components labelled $\sigma$ are
identical. This explains why the construction does not work
for~RB. The crucial lemma is Proposition~\ref{red-subf}(\ref{rs-box}):
If $\Box \phi \in \Gamma(\sigma)$, $R\sigma\sigma'$, and $\sigma'$
occurs in $H'$, then $\phi \in \Gamma(\sigma')$.  Suppose we tried to
define the counter-example $M$ with the the symmetric closure of~$R$ as
its accessibility relation. Then we would have to change the
definition of reduction so as to not only add $\phi$ to the antecedent
of $H(\sigma')$ if $\Box\phi \in \Gamma(\sigma)$ (with
$R\sigma\sigma'$) but also vice versa. Then
Proposition~\ref{sigma-tree}(\ref{st-unchanged}) would no longer hold.
Hence the prospects of extending the method of proving cut-free
completeness to RB are dim.

Recently, a cut-free complete linear nested sequent system for~B has
been developed by \citet{GoreLellmann2019}, though it requires the
introduction of a new modal rule and so does not obey Do\v{s}en's
principle. A cut-free tree hypersequent system for symmetric logics
has also been developed by \citet{Poggiolesi2010}. The structure of
tree hypersequents provides more structural flexibility, and is thus
capable of accommodating symmetric frame properties. However, again we
see Do\v{s}en's principle violated.

More surprisingly, constructing transitive counter-examples for
R4-unprovable hypersequents also causes difficulties. Here the
problem is different and resides in the ``destructive'' nature of the
unrestricted EW rule. Suppose we were going to define a transitive
counter-example $\mathfrak{M}$ using the transitive closure of~$R^1$. Then the
definition of reduction would have to take into account not just
immediate predecessors of $\sigma$ (as the $\Box$L reduction does),
but \emph{any} predecessor of~$\sigma$, i.e., we would define
\begin{align*}
& G \HSeq \Box\phi, \Gamma' \Seq[\sigma'] \Delta' \HSeq G'' \HSeq \phi,
\Gamma \Seq[\sigma] \Delta \HSeq G'
\intertext{to be a 4-reduct of}
& G \HSeq \Box \phi, \Gamma' \Seq[\sigma'] \Delta' \HSeq G'' \HSeq \Gamma
\Seq[\sigma] \Delta \HSeq G'
\intertext{However, if $G''$ is not empty, the unprovability of a
hypersequent does not guarantee the unprovability of its 4-reduct. The best
we can do is guarantee the unprovability of}
& G \HSeq \Box\phi, \Gamma' \Seq[\sigma'] \Delta \HSeq \phi, \Gamma
\Seq[\sigma] \Delta \HSeq G'
\end{align*}
using the EW and $\Box$L rules.  But now the new reduct is no longer
an extension of the original hypersequent, and so
Proposition~\ref{full-red-unprov} fails. The problem, in short, is
that EW destroys information that is required in the subsequent
reduction of a hypersequent and of its successor hypersequents.

The problem can be circumvented by using rules other than EW to deal
with transitivity. One could strengthen the $\Box$L to the rule
\begin{prooftree}
  \def\fCenter{\HSeq}
\Axiom$G \HSeq \phantom{\Box\phi,} \Gamma' \Seq \Delta' \fCenter G'' \HSeq
\phi, \Gamma \Seq \Delta \HSeq G'$
  \RightLabel{$\Box$L4}
\UnaryInf$G \HSeq \Box \phi, \Gamma' \Seq \Delta' \fCenter G'' \HSeq
\phantom{\phi,} \Gamma \Seq \Delta \HSeq G'$
\end{prooftree}
or add a transitivity rule like
\begin{prooftree}
  \def\fCenter{\HSeq}
\Axiom$G \HSeq \phantom{\Box \phi,} \Gamma' \Seq \Delta' \fCenter G'' \HSeq
\Box\phi, \Gamma \Seq \Delta \HSeq G'$
  \RightLabel{Tran}
\UnaryInf$G \HSeq \Box \phi, \Gamma' \Seq \Delta' \fCenter G'' \HSeq
\phantom{\Box\phi,} \Gamma \Seq \Delta \HSeq G'$
\end{prooftree}
In both cases, the unprovability of a hypersequent would guarantee the
unprovability of its reducts which would furthermore be extensions of
them. However, the resulting calculi no longer satisfy Do\v sen's
Principle, since the new rules are not purely (external) structural
rules. It remains an open question whether or not a cut-free complete
hypersequent system for~B or S4 which obeys Do\v{s}en's principle can
be developed. In fact, it is open if RB, R4, and RS4 are cut-free
complete.\footnote{Obvious potential counterexamples are the axioms $B$
and~4, both of which, however, have cut-free proofs (see
Table~\ref{table:cf-proofs}).}

\begin{table}[t]
\caption{Cut-free proofs in RB and R4 of B and 4}
\label{table:cf-proofs}
\def\fCenter{}
\[\begin{array}{@{}ll@{}}
\hline\hline \\
\Axiom$\fCenter \phantom{\lnot A,{}}A \Seq A$
\RightLabel{$\lnot$L}
\UnaryInf$\fCenter \lnot A, A \Seq$
\RightLabel{EWL}
\UnaryInf$\quad \Seq \quad \HSeq \fCenter \lnot A, A \Seq$
\RightLabel{$\Box$L}
\UnaryInf$\Box \lnot A \Seq \quad \HSeq \fCenter \phantom{\lnot A,{}}A \Seq$
\RightLabel{Sym}
\UnaryInf$A \Seq \quad \HSeq \fCenter \,\,\,\Box \lnot A \Seq $
\RightLabel{$\lnot$R}
\UnaryInf$A \Seq \quad \HSeq \fCenter \phantom{\,\,\,\Box \lnot A} \Seq \lnot \Box \lnot A $
\RightLabel{$\Box$R}
\UnaryInf$A \Seq \Box \lnot \Box \fCenter\lnot A $
\DisplayProof
&
\Axiom$\fCenter \phantom{\HSeq} A \Seq A$
\RightLabel{EWL}
\UnaryInf$\quad \Seq \phantom{\Box A} \fCenter \HSeq A \Seq A$
\RightLabel{$\Box$L}
\UnaryInf$\Box A\Seq \phantom{\Box A} \fCenter \HSeq \phantom{A} \Seq A$
\RightLabel{EW}
\UnaryInf$\Box A\Seq \phantom{\Box\Box A} \HSeq \quad \Seq \phantom{\Box A} \fCenter \HSeq \phantom{A} \Seq A$
\RightLabel{$\Box$R}
\UnaryInf$\Box A \Seq \phantom{\Box\Box A} \HSeq \quad \Seq \Box A \fCenter$
\RightLabel{$\Box$R}
\UnaryInf$\Box A \Seq \Box\Box A \phantom{{} \HSeq \quad \Seq \Box A} \fCenter$
\DisplayProof
\\ {} \\
\hline\hline
\end{array}\]
\end{table}

\section{Conclusion}

In this paper we have studied the hypersequent systems of
\citet{Parisi2017}, which extend the hypersequent system for~S5 due to
\citet{Restall2009} to other modal logics. Though these systems
require only two modal rules for~K and its extensions, there are some
issues that arise with this approach to modal hypersequents. In
particular, our method for showing cut-free completeness fails for the
systems RB and RS4. All known cut-free complete hypersequent calculi
for RB and RS4 add additional rules that manipulate modal formulas,
and so violate Do\v{s}en's principle. We have also noted that the
system RS5 is not entirely modular. While it has been shown to be
cut-free complete, in order to obtain modularity by replacing EE with
Sym, we lose the cut-free completeness result. 

\subsection*{Acknowledgements}

The authors would like to thank the reviewers for this paper for their
helpful comments. The results of in this paper are based on
\citet{Burns2018} and were presented at the Society for Exact
Philosophy and the Melbourne Logic Group. We would like to thank the
audiences for their helpful comments and criticisms. A special thanks
to Andrew Parisi, who presented us with an early version of his
foundational work on this topic and gave us the opportunity to expand
upon it.

\bibliographystyle{asl}

\end{document}